\documentclass{elsart}
\usepackage{graphicx,epsfig}
\usepackage{hyperref}
\usepackage{enumerate}

\begin{document}

\begin{frontmatter}
\title{Hadwiger's conjecture for Proper Circular Arc graphs}

\author[IISc]{Naveen Belkale}, 	
\ead{naveenb@csa.iisc.ernet.in,bnaveenbhat@gmail.com}
\author[IISc]{Sunil Chandran} 
\ead{sunil@csa.iisc.ernet.in}
\address[IISc]{CSA Department, Indian Institute of Science, Bangalore, India. Pin-560012}

\begin{abstract}
Circular arc graphs are graphs whose vertices can be represented as arcs on a 
circle such that any two vertices are adjacent if and only if their 
corresponding arcs intersect. Proper circular arc graphs are graphs which have 
a circular arc representation where no arc is completely contained in any other 
arc. Hadwiger's conjecture states that if a graph $G$ has chromatic number $k$, 
then a complete graph on $k$ vertices is a minor of $G$. We prove 
Hadwiger's conjecture for proper circular arc graphs.
\end{abstract}

\end{frontmatter}

\section{Introduction}

Circular arc graphs are graphs whose vertices can be represented as arcs on a 
circle such that any two vertices are adjacent if and only if their 
corresponding arcs intersect. Circular arc graphs form useful mathematical 
objects with many practical applications such as in 
{ \em Genetic research~}\cite{Stahl}, {\em Traffic control}~\cite{Stouffers},
{\em Compiler
design}~\cite{Tuck} and {\em Statistics}~\cite{Hubert}. 

Circular arcs were first discussed in \cite{Klee} as a natural generalization
of interval graphs (defined analogously, but using intervals on a line instead
of arcs on a circle) and they have since been studied extensively 
\cite{Tuck,Tuck1,Tuck2,Garey,Gavril,Hsu,Masuda}. 
See Golumbic~\cite{Golu} for a brief introduction on circular arc graphs.
The {\em circular arc coloring problem} consists of finding a minimum coloring
of a set of arcs of a circle such that no two intersecting arcs have the
same color. Tucker~\cite{Tuck} gave a simple 2-approximation algorithm for 
coloring circular arc graphs and conjectured that $\frac{3}{2} \omega(F)$ 
colors are sufficient to color any family $F$ of
arcs, where $\omega(F)$ represents the size of the maximum set of pairwise 
intersecting arcs. Karapetyan~\cite{Kara} proved Tucker's conjecture. Further 
analysis of Tucker's algorithm by Pabon~\cite{Pabon} gave a tighter upper bound 
for coloring the arcs based on the {\em circular cover} of the family of 
arcs. Circular cover of a family $F$ of arcs is the minimum number of arcs 
required to cover the entire circle.

A graph $G$ is a {\em proper} circular-arc graph if there is a circular-arc 
representation of $G$ such that no arc is properly contained in any other. 
The set of arcs in a proper circular arc representation is called a 
{\em family of proper arcs}.
Because of this additional restriction, some of the difficult problems for 
circular arc graphs in general becomes easier for proper circular arc graphs.
For example, even though coloring of circular arc graphs in general was
shown to be NP hard by  Garey et al.~\cite{Garey}, proper circular arc graphs
can be colored in polynomial time as shown by Orlin, Bonuccelli and 
Bovet~\cite{Orlin}.
Proper circular arc graphs are a strict subclass of circular arc graphs.
See Tucker~\cite{Tuck1,Tuck2} for characterization of proper circular arc graphs.
In this paper, we prove Hadwiger's conjecture for proper circular arc graphs.

\begin{defn}
A \textbf{vertex coloring} of a graph $G=(V,E)$ is a 
map $c: V \rightarrow \{1, \ldots, k\}$ such that $c(v) \neq c(w)$ whenever $v$ 
and $w$ are adjacent. The smallest integer $k$ such that $G$ has a vertex 
coloring $c:V \rightarrow \{ 1, \ldots, k\}$ is called the \textbf{chromatic 
number} of $G$. Chromatic number of a graph $G$ is denoted by $\chi(G)$.
\end{defn}
\begin{defn}
\textbf{Contraction of an edge $e=(x,y)$} is the replacement
of vertices $x$ and $y$ with a new vertex $z$ whose incident edges are 
exactly  those edges
other than $e$ that were incident to at least one of  $x$ or $y$.
\end{defn}
\begin{defn}
A \textbf{Minor} of a graph $G$ is obtained by a sequence of applications 
of the following three elementary operations on graph $G$ 
\begin{enumerate}
\item Deletion of a vertex
\item Deletion of an edge
\item Contraction of an edge 
\end{enumerate}
\end{defn}
\begin{defn}
\textbf{Complete graph} of order $n$ is a graph on $n$ vertices where all the
vertices are pairwise adjacent and is denoted by $K_n$.
\end{defn}

See Diestel~\cite{Diest} for more information on minors.

In 1943, Hadwiger made the famous conjecture linking the chromatic number of a
graph with its clique minor:
\begin{conj}[Hadwiger~\cite{Hadwiger43}]
If a graph $G$ has chromatic number $\chi(G)=r$, then $K_r$ is a minor of $G$.
\end{conj}

The conjecture is easy to see for $r=1,2,3$. The case $r=4$ was proved by 
Dirac~\cite{Dirac3}. The Hadwiger's conjecture for any fixed $r$ is equivalent
to the assertion that every graph without a $K_r$ minor
has a $(r-1)$ coloring. Therefore, the case $r=5$ implies four color theorem
because any planar graph has no $K_5$ minor. On the other hand, Hadwiger's 
conjecture for the case $r=5$ follows from the four color theorem and a 
structure theorem of Wagner~\cite{Wagner2}.

For $r=6$, Dirac~\cite{Dirac4} proved that if the chromatic number of a graph
$G$ is 6, then $G$ can be contracted into ${K_6}^-$, a complete graph on 6 nodes
with one edge missing. Robertson, Seymour and Thomas~\cite{RobertsonSeymourThomas93} have obtained
a proof for $r=6$. They showed that every minimum contraction-critical graph
(a graph is said to be contraction-critical if for every proper minor $H$ of 
$G$, $\chi(H) < \chi(G)$), different from the complete graph $K_6$, is an
{\em apex} graph which has a special vertex $v$ the removal of which results in
a planar graph. As a result, Hadwiger's conjecture for $r=6$ reduces to the 
four color theorem. The case for $r=7$ onwards is still open and the best known
results for $r=7,8,9$ are due to Jakobsen~\cite{Jack1,Jack2}. He proved that a 
$k$-chromatic graph can be contracted to ${K_7}^{--}$, ${K_7}^-$ and $K_7$ 
respectively for $k=7,8$ and $9$.  Recently, Kawarabayashi and 
Toft~\cite{Toft2} proved that any $7$-chromatic graph has $K_7$ or $K_{4,4}$ as 
a minor.

Since Hadwiger's conjecture in the  general case  seems to be too difficult, 
it is interesting to prove it  for special classes of graphs.
Reed and Seymour~\cite{ReedSeymour} studied the Hadwiger's conjecture in the
case of {\em line graphs}. A line graph $L(G)$ of $G$ is the graph on the edge
set $E(G)$ in which $x,y \in E(G)$ are adjacent as vertices in $L(G)$,
 if and only if they are adjacent as edges in $G$. They showed that for every 
loop less graph $G$(possibly with parallel edges), its line graph $L(G)$ 
satisfies the Hadwiger's conjecture. For a detailed history of Hadwiger's 
conjecture as well as an account of recent developments in that area the reader 
is referred to the survey by Toft~\cite{Toft}.

A graph is called perfect if for every induced subgraph $H$ of $G$, 
$\chi(H)=\omega(H)$ where $\omega(H)$ is the order (i.e.,
 the number of vertices) of the largest complete 
subgraph of $H$. $\omega(H)$ is  also known as the   {\em clique number}
 of $H$. Interval graphs are 
{\em perfect} graphs. See Golumbic~\cite{Golu} for more information on 
interval graphs.  Hence Hadwiger's conjecture trivially holds true for interval
graphs. For the class of Circular arc graphs, which is a generalization of 
interval graphs, Hadwiger's conjecture is still open. Therefore, it is natural 
to study Hadwiger's conjecture for circular arc graphs. In this paper, we
prove Hadwiger's conjecture for proper circular arc graphs.

\begin{rem} {Hajos conjecture for proper circular arc graphs:}
\end{rem}
It may be a matter of curiosity to note that the class of proper circular arc
graphs doesn't satisfy the closely related Hajos conjecture. This is because
the counterexample to Hajos conjecture given by Catlin~\cite{Catlin} is a proper circular 
arc graph. This gives us one more reason to verify Hadwiger's conjecture
for proper circular arc graphs.

\section{Preliminaries}

For a circular arc graph $G$, without loss of generality we can assume 
that the family of arcs $F$ representing $G$ has all its arc endpoints 
distinct. Also, without loss
of generality we can assume that no arc in the circular arc 
representation of a circular arc graph spans the whole circle
Similarly, we can assume that no arc is formed of a single point. 
These assumptions
also hold true for proper circular arc graphs.

\begin{defn}
Given a family $F$ of arcs, the \textbf{overlap set} of point $p$
 on the circle 
is the set of all arcs that contain the given point $p$ and is denoted by  
$\mathcal{O}(p)$. The overlap set with the largest number of arcs in it is 
called a
\textbf{maximum overlap set} and its cardinality is denoted by $r_{sup}$. An
overlap set with the smallest number of arcs in it is called a \textbf{minimum
overlap set} and its cardinality is denoted by $r_{inf}$.
\end{defn}

It is easy to see that the arcs in an overlap set induce a clique in the 
corresponding circular arc graph. From now on, when there is no ambiguity we use
the term ``arcs'' and ``vertices'' interchangeably. For example, we use the same
labels to refer to the vertices in a circular arc graph as well as the 
corresponding arcs in its circular arc representation that is being considered.

\begin{defn}
Each arc in $F$ has two endpoints. The \textbf{left endpoint} $l(u)$ 
$($respectively \textbf{right endpoint $r(u)$}$)$  of arc $u$ is the first 
endpoint of $u$ encountered in an {\em anticlockwise} $($respectively
{\em clockwise}$)$ traversal from any interior point of $u$. 
$($ Recall that the circle itself is not considered as an arc. A single
point is also not considered as an arc. Thus the definition makes sense,
and each arc has a left end point as well as a right end point.$)$
\end{defn}

\begin{defn}
An arc $v$ is said to be \textbf{clockwise adjacent} to an arc $u$ if 
$v$ belongs to the overlap set of $r(u)$, i.e. $v \in \mathcal{O}(r(u))$. An 
arc $v$ is said to be \textbf{anticlockwise adjacent} to an arc $u$ if $v$ 
belongs to the overlap set of $l(u)$, i.e. $v \in \mathcal{O}(l(u))$.
\end{defn}

For a family of arcs, an arc $u$ can be adjacent to an arc $v$ without being
{\em clockwise adjacent} or {\em anticlockwise adjacent} to $v$ by being 
properly contained in $v$.  For a family of {\em proper} arcs, if an arc $u$ 
intersects arc $v$ then it must be either {\em clockwise} or 
{\em anticlockwise adjacent} to $v$.\\
\begin{note} In a family of proper arcs, if an arc $u$ is clockwise adjacent to 
an arc $v$, then arc $v$ is anticlockwise adjacent to arc $u$.
\end{note}

\begin{lem} \label{lem:lem24}
An arc in a family of proper arcs has at most $r_{sup}-1$ arcs and at least 
$r_{inf}$ arcs clockwise adjacent to it. Similarly, an arc has at most 
$r_{sup}-1$ and at least $r_{inf}$ arcs anticlockwise adjacent to it.
\end{lem}
\begin{pf}
There can be at most $r_{sup}-1$ arcs that are clockwise adjacent to an arc $u$
since cardinality of overlap set at the right endpoint of $u$, 
$|\mathcal{O}(r(u))| \leq r_{sup}$ and $u$ is also part of this overlap set.
Similarly, we can prove for the anticlockwise direction using $l(u)$.

For an arc $u$, consider the point $p$ just after right endpoint $r(u)$ of arc 
$u$ in clockwise direction.
As the cardinality of the minimum overlap set is $r_{inf}$, we have 
$|\mathcal{O}(p)| \geq r_{inf}$. Clearly, $\mathcal{O}(p) \subset 
\mathcal{O}(r(u))$ and hence each arc in $\mathcal{O}(p)$ is clockwise adjacent
to $u$. Therefore, there are at least $r_{inf}$ arcs that are clockwise
adjacent to $u$. Similarly we can argue for the anticlockwise direction by
looking at a point $p$ just after $l(u)$ in anticlockwise direction. \qed
\end{pf}

The minimum number of colors needed to color a family $F$ of arcs such that
no two intersecting arcs have the same color is its {\em chromatic number} 
$\chi(F)$. As the arcs correspond to vertices in the corresponding 
circular arc graph $G$, we have $\chi(F)=\chi(G)$. A straightforward upper bound
on the chromatic number of $F$ is
\begin{lem}[Tucker~\cite{Tuck}] \label{lem:lem25}
For a family $F$ of arcs, $\chi(F) \leq r_{sup} + r_{inf}$.
\end{lem}

This can be easily seen from the fact that removing all the arcs in any
minimum overlap set will result in a family of arcs  that correspond to an
interval graph which can then be colored using at most $r_{sup}$ colors. 
(Recall
that interval graph is a perfect graph and the cardinality of maximum overlap 
set corresponds to the clique number in an interval graph.)

\begin{defn} \label{def:def26}
For a family $F$ of arcs, \textbf{circular cover $l(F)$} is the smallest 
cardinality  of
any subset of arcs of $F$ needed to cover the circle.
\end{defn}

Note that circular cover is defined for a family of arcs only if a finite 
number of arcs in the family can cover the circle. If no number of arcs can 
cover the entire circle then the corresponding graph is an interval graph and
we know the chromatic number equals clique number for interval graphs
and thus Hadwiger's conjecture trivially holds. Hence in the following 
sections we will assume that the circular cover is  defined and finite.

\begin{thm}[Tucker~\cite{Tuck}] \label{thm:thm26}
If the circular cover of a family $F$ of arcs $l(F) \geq 4$, then $\chi(F) \leq 
\frac{3}{2} r_{sup}$
\end{thm}

\section{Hadwiger's conjecture for proper circular arc graphs}

Suppose, for contradiction, assume that the class of proper circular
arc graphs does not satisfy Hadwiger's conjecture. Then, consider a proper
circular arc graph on the minimum possible number of vertices that 
is a counterexample to Hadwiger's conjecture.

\begin{table} 
\begin{tabular}{|p{0.08\linewidth}p{0.85\linewidth}|} 
\hline
$\mathcal{G}$ & a proper circular arc graph on the minimum possible  number of 
vertices  that is a counterexample to Hadwiger's conjecture\\
$\chi(\mathcal{G})$ & chromatic number of the graph $\mathcal{G}$\\
$\delta(\mathcal{G})$ & minimum degree of the graph $\mathcal{G}$\\
$r$ & number of arcs in maximum overlap set of family of proper arcs $\mathcal{G}$\\
$O$ & a maximum overlap set of family of proper arcs $\mathcal{G}$\\
$n$ & number of vertices in the graph $\mathcal{G}$\\
$x$ & integer such that $\chi(\mathcal{G}) = r + x$ \\
$k$ & integer such that $n = r + k$\\
Note:& Notation $\mathcal{G}$ is used both for the proper circular arc graph as well as for the proper circular arc representation under consideration for the graph \\
\hline
\end{tabular}
\caption{Notation for the minimum counterexample}
\label{table:table1} 
\end{table}

We will use the notations given in Table~\ref{table:table1} throughout the 
remaining part of this section. 
If $\chi(\mathcal{G})=r$, then Hadwiger's conjecture trivially holds
true for $\mathcal{G}$ as all the arcs of the maximum overlap set form a 
clique.  Therefore, we have $x >0$.  Also, if $n=r$, then $\mathcal{G}$ will be 
a complete graph and hence not a counterexample to Hadwiger's conjecture.
Hence, $k>0$. We fix a proper circular arc representation for graph 
$\mathcal{G}$ which will also be referred using the same notation $\mathcal{G}$.

The following theorem is well known in the literature 
regarding Hadwiger's conjecture (See Kotlov~\cite{Kotlov} for an alternate 
proof.)
\begin{thm}[Gallai~\cite{Galla1}] \label{thm:thm31}
If $G$ on $n$ vertices is the only counterexample to Hadwiger's conjecture 
among its induced subgraphs then $\chi(G) \leq \lceil n/2 \rceil$.
\end{thm}

Any induced subgraph of a proper circular arc graph is also a proper circular
arc graph. This is because removing a vertex in the graph is equivalent to 
removing an arc in the corresponding circular arc representation.  As 
$\mathcal{G}$ is a
proper circular arc graph on the smallest possible number of vertices that is a 
counterexample to Hadwiger's conjecture, it is the only counterexample 
to Hadwiger's conjecture among its induced subgraphs. Therefore by 
Theorem~\ref{thm:thm31},

\begin{lem} \label{lem:lem32}
$\chi(\mathcal{G})\leq \lceil n/2 \rceil$
\end{lem}

A graph $G(V,E)$ is said to be color critical if $\chi(G-v) < \chi(G)$ for every
$v \in V$. Every graph has an induced subgraph that is color critical. 
For a color critical graph $G$, we have the minimum degree of the graph, 
$\delta(G) \geq \chi(G)-1$(See West~\cite{West} for the proof).
It is easy to see that the minimum counterexample to Hadwiger's conjecture 
$\mathcal{G}$ must be  {\em color critical}. Therefore we have,
\begin{lem} \label{lem:lem33}
$\delta(\mathcal{G})\geq r+x-1$
\end{lem}

\begin{lem} \label{lem:lem2x}
$2x \leq r$
\end{lem}
\begin{pf}
We will first show that $\chi(\mathcal{G}) \leq \frac{3}{2}r$.
If the circular cover (see definition~\ref{def:def26}) $l(\mathcal{G}) \geq 4$, 
then by Theorem~\ref{thm:thm26}, we have $\chi(\mathcal{G}) \leq \frac{3}{2}r$. 
In a family of proper arcs, the set of  arcs that are clockwise
adjacent to a given arc is exactly the set of arcs that have their 
left endpoints inside that  arc. Therefore,
by Lemma~\ref{lem:lem24}, every 
arc in $\mathcal{G}$ can have left endpoints of at most $r-1$ arcs in it.
If $l(\mathcal{G}) \leq 3$, then there exists three arcs say $x$, $y$ and $z$ 
such that the union of these arcs cover the entire circle. Since each arc 
including $x$,$y$ and $z$ should have its left endpoint in the interior of 
at least one of these three
arcs, it follows that $n \leq 3(r - 1)$. As $\mathcal{G}$ is a minimum 
counterexample by Lemma~\ref{lem:lem32}, we have $\chi(\mathcal{G}) \leq 
\lceil \frac{3r-3}{2} \rceil \leq \frac{3}{2}r$. Hence even if $l(\mathcal{G}) \leq 3$, we have $\chi(\mathcal{G}) \leq \frac{3}{2}r$. Recalling that 
$\chi(\mathcal{G})=r+x$, we get $2x \leq r$. \qed
\end{pf}

\begin{lem} \label{lem:lemr2x}
$k \geq r+2x-1$
\end{lem}
\begin{pf}
From Lemma~\ref{lem:lem32} we have $r + x \leq \lceil n/2 \rceil \leq (n+1)/2
= (r+k+1)/2$  and the lemma follows. \qed
\end{pf}
\begin{cor} \label{cor:cor34}
$k \geq 4x-1$
\end{cor}

\begin{lem} \label{lem:fr}
For a family $F$ of proper arcs, if we are traversing in the clockwise 
direction from a point $p$ on the circle, the right endpoints of all the arcs 
in the overlap set of $p$ would be encountered before the right endpoints of 
any other arc in $F$.
\end{lem}
\begin{pf}
When we traverse along the circle in clockwise direction from a point $p$, if 
we encounter the right endpoint of an arc $u$ before encountering its left 
endpoint, then $u$ is in the overlap set of $p$, i.e. $u \in \mathcal{O}(p)$.
Instead, as we traverse along the circle in clockwise direction from $p$, if
we encounter the left endpoint of arc $v$  before its right endpoint, then $v 
\notin \mathcal{O}(p)$. If the right endpoint of $v$ ($v \notin 
\mathcal{O}(p)$) also occurs before the right endpoint of at least one arc $u$ 
($u \in \mathcal{O}(p)$), then  moving in clockwise direction from $l(u)$ 
we encounter the point $p$ and then both the left and right endpoints of $v$ 
before we encounter $r(u)$ which implies
that $v$ is entirely contained in arc $u$ which is in contradiction to proper 
circular arc property. Hence the lemma follows. \qed
\end{pf}

\begin{cor} \label{cor:fr1}
(a) If an arc $u$ is clockwise adjacent to an arc $v$, then all the arcs whose
right endpoints are encountered after the right endpoint of $v$ and before
the right endpoint of $u$ in clockwise direction are also clockwise adjacent
to arc $v$. 

(b) Similarly, if an arc $u$ is anticlockwise adjacent to an arc $v$,
then all the arcs whose right endpoints are encountered after the right endpoint
of $u$ and before the right endpoint of $v$ in clockwise direction are also
anticlockwise adjacent to arc $v$.
\end{cor}

\begin{figure}
\centering
\includegraphics[width=180pt]{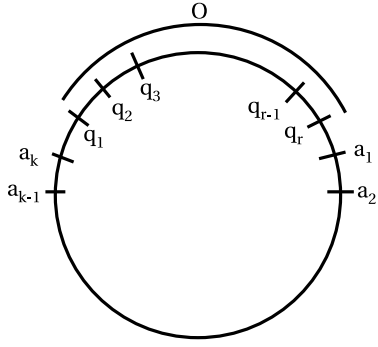}
\caption{Labeling scheme for arcs in $\mathcal{G}$}
\label{fig:fig1}
\end{figure}

In the circular arc representation of $\mathcal{G}$, identify a point $p$ such 
that $\mathcal{O}(p)=O$, the maximum overlap set. Traverse the circle in the 
clockwise direction
starting from the point $p$ labeling the arcs in the order in which their right
endpoints are encountered. Note that the first $r$ arcs to be labeled
are from the maximum overlap set $O$ by Lemma~\ref{lem:fr}. Let the first $r$ 
arcs(i.e. the arcs in $O$) be labeled as $q_1, q_2, \ldots, q_r$ and the 
remaining $k$ arcs be labeled $a_1, a_2, \ldots, a_k$. See Fig~\ref{fig:fig1} 
for the labeling scheme.  By Corollary~\ref{cor:cor34}, we have $k \geq 4x-1$ 
and therefore $\{ a_1, a_2, \ldots, a_x \} \cap \{a_{k-x+1}, \ldots, a_k\} = 
\emptyset$. 

Now based on Corollary~\ref{cor:fr1} 
 we can infer more about the adjacency 
relationships between arcs in $\mathcal{G}$.
\begin{cor} \label{cor:fr2}
In the labeling scheme defined above,
\begin{enumerate}
\item If an arc $a_j$ is clockwise adjacent to an arc $q_i$, then there are at 
least $r-i+j$ arcs, namely $\{q_{i+1}, \ldots, q_r, a_1, \ldots, a_j\}$ 
clockwise adjacent to arc $q_i$.
\item If an arc $a_j$ is clockwise adjacent to an arc $a_i$, then $j>i$ and 
there are at least $j-i$ arcs, namely $\{a_{i+1}, \ldots, a_j\}$ clockwise 
adjacent to arc $a_i$.
\item If an arc $a_j$ is anticlockwise adjacent to an arc $a_i$, then $j<i$ and
there are at least $i-j$ arcs, namely $\{a_j, \ldots, a_{i-1}\}$ anticlockwise 
adjacent to arc $a_i$.
\item If an arc $a_j$ is anticlockwise adjacent to an arc $q_i$, then there 
are at least $k+1-j+i-1=k-j+i$ arcs, namely $\{a_j, \ldots, a_k, q_1, \ldots, 
q_{i-1}\}$ anticlockwise adjacent to arc $q_i$.
\end{enumerate}
\end{cor}

\begin{lem} \label{lem:lem35}
Arc $a_{i+j}$ must be clockwise adjacent to arc $a_i$ for $1 \leq j \leq x$ 
and $1 \leq i \leq k-x$.
\end{lem}
\begin{pf}
As $\mathcal{G}$ is a family of proper circular arcs, any arc $u$ in 
$\mathcal{G}$ has at least $r_{inf}$ arcs clockwise adjacent to it by 
Lemma~\ref{lem:lem24}. In other words, $|\mathcal O (r(u))| \ge r_{inf}$. 
Therefore by Lemma~\ref{lem:fr}, the $r_{inf}$ arcs 
whose right endpoints are encountered after the right endpoint of $u$ in
clockwise direction are 
in the overlap set of $r(u)$ and hence clockwise adjacent to $u$.
By Lemma~\ref{lem:lem25} and the fact that $\chi(\mathcal{G})=r+x$ we have
$x \leq r_{inf}$. Therefore, the $x$ arcs whose right endpoints are encountered
after the right endpoint of an arc $u$ in clockwise direction are
clockwise adjacent to $u$. By the labeling scheme described above, for any 
arc $a_i$ where $1 \leq i \leq k-x$, these $x$ arcs whose right endpoints
are encountered after the right endpoint of $a_i$ would be labeled $a_{i+1}, \ldots, a_{i+x}$. Therefore $a_{i+1}, \ldots, a_{i+x}$ are clockwise 
adjacent to $a_i$. \qed
\end{pf}

We now define a {\em good path set} with respect to the circular arc 
representation $\mathcal{G}$ and the labeling scheme described above.
\begin{defn} 
A \textbf{good path set} is a set of $x$ vertex disjoint paths $P_1, P_2, \ldots
, P_x$ such that each $P_i$ starts at $a_i$ and ends at $a_{k-x+i}$ and 
$P_i \cap O = \emptyset$ where $O$ is the maximum overlap set.
\end{defn}

\begin{lem} \label{lem:lem37}
$\mathcal{G}$ does not contain a good path set.
\end{lem}
\begin{pf}
If $\mathcal{G}$ contains a good path set then we will show a 
$r+x$ (recall that $\chi(\mathcal{G})=r+x$)
clique minor leading to a contradiction. For this, we will first show that 
every arc in $O$ is adjacent to either $a_i$ or $a_{k-x+i}$ for 
$1\leq i \leq x$.

Arcs in $O$ are labeled $q_1, q_2, \ldots, q_r$.
Suppose we have an arc $q_j$ and an integer $i$, where $1\leq i \leq x$ 
such that $q_j$ is adjacent to neither $a_i$ nor $a_{k-x+i}$. 
What is the degree of $q_j$? Clearly $q_j$ intersects with all
the remaining $r-1$ arcs in $O$, at most $i-1$ arcs from $\mathcal{G}-O$ are
clockwise adjacent to $q_j$, namely $a_1, \ldots, a_{i-1}$  and at most $x-i$ arcs from $\mathcal{G}-O$ are anticlockwise adjacent to $q_i$, namely 
$a_{k-x+i+1},\ldots,a_k$.  Therefore, $degree(q_j) \leq (r-1)+(i-1)+(x-i) = 
r+x-2$ which 
contradicts Lemma~\ref{lem:lem33} by which $\delta(\mathcal{G}) \geq r+x-1$. 
Hence, every arc in $O$ must be adjacent to either $a_i$ or $a_{k-x+i}$ for 
$1 \leq i \leq x$.

Now if we contract each vertex disjoint path $P_i = (a_i, \ldots, a_{k-x+i})$ 
to a single vertex, then such a vertex would be adjacent to all the arcs in 
$O$. Also each of the contracted vertices would be adjacent to each other as 
$a_2, \ldots, a_x$ are clockwise adjacent to $a_1$ by Lemma~\ref{lem:lem35}.  
Also, all the arcs in $O$ are pairwise adjacent to each other as they all belong
to the maximum overlap set. Hence, we have an $r+x$ clique minor if 
$\mathcal{G}$ has a good path set contradicting the assumption that 
$\mathcal{G}$ is a counterexample to Hadwiger's conjecture. \qed
\end{pf}

\begin{lem} \label{lem:lemk}
$k$ is not divisible by $x$
\end{lem}
\begin{pf}
If $k$ is a multiple of $x$ then we show a good path set which contradicts 
Lemma~\ref{lem:lem37}. We define $x$ vertex disjoint paths where each path $P_j$
for $1\leq j\leq x$ is of the form $P_j=(a_j, a_{x+j}, \ldots, a_{(t-1)x+j})$
and $t=k/x$. By Lemma~\ref{lem:lem35}, $a_{l+x}$ must be clockwise adjacent 
to $a_l$ for $1\leq l\leq k-x$ and hence each of $P_j$ is a path.
It is easy to see that each vertex $a_l$ for $1\leq l \leq k$ belongs to
a unique path, namely  $P_j$
where $j=(l-1)(\bmod x)+1$. Moreover, vertex $a_{(t-1)x + j} = a_{k-x+j}$ 
is in $P_j$. Thus, we have a good path set contradicting Lemma~\ref{lem:lem37}. \qed
\end{pf}

\begin{lem} \label{lem:lem38}
$a_i$ is not clockwise adjacent to $a_{i-2x+1}$ in $\mathcal{G}$ for any 
$2x \leq i \leq k-x$
\end{lem}
\begin{pf}
From Lemma~\ref{lem:lemk}, $k$ is not a multiple of $x$. Let $t = \lfloor
k/x \rfloor$ and $b = k-tx$. Suppose we have an $a_i$ which is clockwise 
adjacent to $a_{i-2x+1}$ where $2x \leq i \leq k-x$. We will demonstrate a good 
path set which contradicts Lemma~\ref{lem:lem37}.

Define a successor function $s: \{a_i : 1 \le i \le k - x\} \rightarrow
  \{ a_j : x < j \le k \}$ as follows: 
\begin{eqnarray}
s(a_l) &=& a_{l+x} \textrm{ for } 1 \leq l \leq i-2x \nonumber \\
s(a_l) &=& a_{l+x+b} \textrm{ for } i-2x < l \leq i-x-b \nonumber \\
s(a_l) &=& a_{l+b} \textrm{ for } i-x-b < l \leq i-x \nonumber \\
s(a_l) &=& a_{l+x} \textrm{ for } i-x < l \leq k-x \nonumber
\end{eqnarray}

(The reader may note that, when $i = 2x$, the range 
defined by $1 \le l \le i - 2x$ is empty. Thus when $i=2x$, the
function $s$ is completely defined by the last three equations given above.)

\begin{figure}
\centering
\includegraphics[width=370pt]{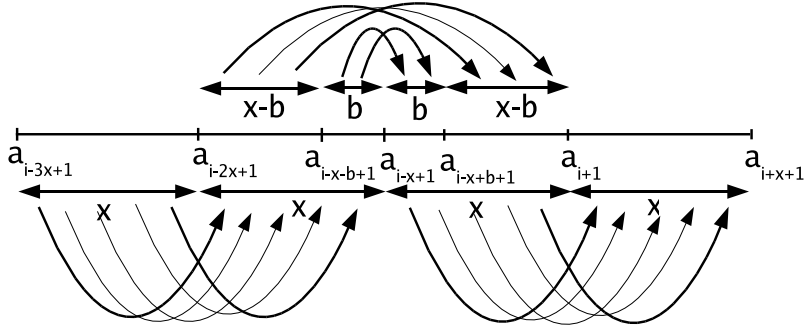}
\caption{Diagram showing successor vertices using arrow}
\label{fig:fig2}
\end{figure}

A {\em successor vertex} of a vertex is the vertex obtained by applying
the successor function on the given vertex.
See Fig~\ref{fig:fig2} for a demonstration of the successor function where
the arrow points from a vertex to its successor vertex. 
We will now demonstrate $x$ vertex disjoint paths of a good path set 
using the successor function defined above, contradicting Lemma 
\ref {lem:lem37}. The required result would follow immediately. 

Define $x$ sequences of vertices $P_1,P_2,\ldots,P_x$ as follows:
Each sequence $P_j$ for $1 \leq j
\leq x$ begins at $a_j$. Any other vertex in the sequence is
determined by applying the successor function on the previous vertex in the 
sequence. For example, $P_j = (a_j, s(a_j), s(s(a_j)), \ldots)$.
The sequence ends when we encounter a vertex for which the successor
function is not defined. Note that the successor function is defined
for every vertex $a_h$ in the range $1 \le h \le k-x$. Moreover, by
applying the successor function on a vertex $a_h$ (where
$1 \le h \le k-x$)  we always get another vertex $a_k$ whee $k > h$.
Therefore the sequences defined above are of  finite length, 
and the last vertex of each of these sequences belong to the set 
$\{a_{k-x+1}, \ldots, a_k\}$,  for which the  successor function is not 
defined.

In order to  demonstrate  that the sequences $P_1,P_2,\ldots,P_x$
indeed form  a good path set, we prove the following claims: \\
\textbf{claim 1.} $s(a_l)$ is clockwise adjacent to $a_l$ for 
$1 \leq l \leq k-x$. Thus for $1 \le j \le x$, the sequence 
$P_j$ forms path.\\
\textbf{proof.} 
\begin{enumerate}
\item This is easy to see in the range $1 \leq l \leq i-2x$ and in 
the range $i-x-b < l \leq k-x$ by 
Lemma~\ref{lem:lem35}. 
\item In the range $i-2x < l \leq i-x-b$, we have $s(a_l)= a_{l+x+b}$. 
If $a_{i-2x+1}$ is anticlockwise
adjacent to $a_i$ then since $i-2x+1 \leq l < i$, $a_l$ is also 
anticlockwise adjacent to $a_i$ by Corollary~\ref{cor:fr2}.3. 
We have $l < l+x+b \leq (i-x-b)+(x+b) \leq i$. Therefore, by 
Corollary~\ref{cor:fr2}.2, as $a_i$ is clockwise adjacent to $a_l$, $a_{l+x+b}$ 
is also clockwise adjacent to $a_l$. Thus in the range
$i-2x < l \leq i-x-b$ we have $s(a_l)$ clockwise adjacent to $a_l$.
\end{enumerate}
Therefore, in the whole range $1 \leq l \leq k-x$, $s(a_l)$ is clockwise
adjacent to $a_l$.

\textbf{claim 2.} The paths $P_1,P_2,\ldots,P_x$ are vertex disjoint\\
\textbf{proof.} Every vertex $a_l$ in the range $1 \leq l \leq k$ 
belongs to exactly
one path. We can easily determine the unique path to which the vertex belongs by
studying the successor function and it is given below
\begin{enumerate}
\item In the range $1 \leq l \leq i-x$ we have $a_l$ in path $P_j$ where 
$j=(l-1)(\bmod x)+1$. This is because if $a_l$ is in $P_j$ for $1\leq l \leq 
i-2x$ then $s(a_l)=a_{l+x}$ is also in $P_j$ and $a_j$ is in $P_j$, $1\leq 
j\leq x$.  For example, vertices $\{a_1,a_{x+1}, a_{2x+1}, \ldots, a_{hx+1}\}$
where $hx + 1 \le i -x$, belong to path $P_1$.

\item 
In the range $i-x < l \leq i-x+b$, the vertex $a_l$ is a successor
of the vertex $a_{l-b}$. Hence it belongs to the same path to which $a_{l-b}$
belongs. As $l-b \leq i-x$, $a_{l-b}$ belongs to $P_j$ where 
$j=(l-b-1)(\bmod x)+1$.

\item 
In the range $i-x+b < l \leq i$, the vertex $a_l$ is a successor of
the vertex $a_{l-x-b}$ and hence belongs to the same path to which $a_{l-x-b}$
belongs.As $l-x-b \leq i-x-b \leq i-x$, $a_{l-x-b}$ belongs to path $P_j$ where 
$j=(l-x-b-1)(\bmod x)+1=(l-b-1)(\bmod x)+1$.

\begin{note}
The reader may observe that the $x$ vertices
$a_{i-2x+1}, \ldots, a_{i-x}$ have their successor vertices rearranged in the
consecutive $x$ vertices $a_{i-x+1}, \ldots, a_i$. See Fig~\ref{fig:fig2} for
the rearrangement of successor vertices in the above range.
\end{note}

\item From (2) and (3), it is clear that $a_l$ for $i-x < l \leq i$ belongs
to $P_j$ where $j=(l-b-1)(\bmod x)+1$.
In the range $i < l \leq k$, $a_l$ is in the same path to which $a_{l-x}$
belongs and $a_{l-x}$ belongs to $P_j$ where $j=(l-x-b-1)(\bmod x)+1 =
(l-b-1)(\bmod x)+1$.

\end{enumerate}

\textbf {Claim 3.} Every path $P_l$, for $1 \le l \le x$ ends
 at $a_{k-x+l}$.\\
\textbf {proof.} For $l$ in the range $1 \leq l \leq x$, vertex $a_{k-x+l}$ 
belongs
to path $P_j$ where $j=(k-x+l-b-1)(\bmod x)+1 = (tx+b-x+l-b-1)(\bmod x)+1=
(l-1)(\bmod x)+1 = l$.  Then, clearly $a_{k-x+l}$ should be the last
vertex of $P_l$ since the successor function is not defined for $a_{k-x+l}$.

From Claims 1,2 and 3, it is clear that $P_1,P_2,\ldots,P_x$ form a 
good path set. 
This contradicts 
Lemma~\ref{lem:lem37} and hence $a_i$ is not clockwise adjacent to 
$a_{i-2x+1}$ in $\mathcal{G}$ for $i$ in the range $2x \leq i \leq k-x$. \qed
\end{pf}

\begin{lem} \label{lem:lem39}
$2x + (k \bmod x) > r$
\end{lem}
\begin{pf}
Let $b= k \bmod x$. By Lemma~\ref{lem:lemk}, we have $1\leq b \leq x-1$.
If $2x + b \leq r$, we will demonstrate a good path set which contradicts 
Lemma~\ref{lem:lem37}.
Let us consider  $x$ vertex disjoint paths 
$P_1, \ldots, P_x$ where $P_j = (a_j, a_{x+j},
\- \ldots, a_{k-2x-b+j})$. Note that $P_j$ is a path because 
$a_{l+x}$ is clockwise 
adjacent to $a_l$ for $1\leq l \leq k-x$ by Lemma~\ref{lem:lem35}. It is 
easy to see that, for $1 \le l \le k-x-b$  the vertex $a_l$ is in $P_j$ if $(l-1)(\bmod x)+1=j$ and in particular,  
$a_{k-2x-b+j}$ is in $P_j$ as $(k-2x-b+j-1)(\bmod x)+1=j$. Now, we will 
show that, for $1 \le j \le x$,  $a_{k-2x-b+j}$ is anticlockwise adjacent 
to $a_{k-x+j}$. Then we get a
good path set by attaching the vertex $a_{k-x+j}$ to the path $P_j$ after 
$a_{k-2x-b+j}$.

From Lemma~\ref{lem:lem38}, $a_{k-2x-b+j}$ is not clockwise adjacent to
$a_{k-4x-b+j+1}$ (we can apply Lemma~\ref{lem:lem38} to $a_{k-2x-b+j}$ because 
from Lemma~\ref{lem:lemr2x}, $k\geq r+2x-1$ and as 
$2x+b\leq r$ it is clear that $k-2x-b+1 \geq k-r+1 \geq 2x$).
Hence, there are at most $2x-2$ arcs that are anticlockwise
adjacent to $a_{k-2x-b+j}$. It follows from Lemma~\ref{lem:lem33} that 
there are at least $(r+x-1) - (2x-2) = r-x+1$ arcs clockwise adjacent to 
$a_{k-2x-b+j}$. 
If $a_{k-x+j}$ is not clockwise adjacent to $a_{k-2x-b+j}$ then there are at 
most $(k-x+j) - (k-2x-b+j+1) = x+b-1 \leq r-x$ (the last inequality following 
from the assumption $2x+b \leq r$) arcs that are clockwise adjacent to 
$a_{k-2x-b+j}$
which is a contradiction. Hence $a_{k-x+j}$ is clockwise adjacent to 
$a_{k-2x-b+j}$ if $2x+b \leq r$ resulting in a good path set as explained above.\qed
\end{pf}

\begin{thm}
Hadwiger's conjecture is true for proper circular arc graphs.
\end{thm}
\begin{pf}
Suppose Hadwiger's conjecture is false for proper circular arc graphs. Then
let $\mathcal{G}$ be a proper circular arc graph on the smallest number of 
vertices that is a counterexample to Hadwiger's conjecture.
We continue to use the notations provided in Table~\ref{table:table1} and the 
labeling scheme shown in Fig~\ref{fig:fig1} for the proof. Since $\mathcal{G}$ 
is a minimum
counterexample to Hadwiger's conjecture it satisfies Lemma~\ref{lem:lem37}, 
Lemma~\ref{lem:lemk}, Lemma~\ref{lem:lem38} and Lemma~\ref{lem:lem39}. We will 
show that these properties will allow us to color $G$ using $r+x-1$ colors 
leading to a contradiction.

Let $t=\lfloor k/x \rfloor$. Let $b = k-tx$. As $k$ cannot be a multiple of $x$
(by Lemma~\ref{lem:lemk}), we have $1 \leq b \leq x-1$. Now we demonstrate a 
vertex coloring of $\mathcal{G}$, $f: V \rightarrow \{1, 2, \ldots, 
r+x-1 \}$. We need slightly different strategies for coloring depending on 
whether $t$ is even or $t$ is odd.

\begin {note}
\label {note:indset}
After proposing a vertex coloring $f: V \rightarrow \{1,2,\ldots,r+x-1\}$,
to prove that it is a valid vertex coloring, our strategy would be to
show that for each $i$, $1 \le i \le r+x-1$, the set of vertices that
are colored $i$ form an independent set in $\mathcal {G}$. 
To make the discussion easy, we first observe the following simple fact:
Let $X=\{q_i,a_{i_1},a_{i_2},\ldots,a_{i_p}\}$ be a set of vertices of
$\mathcal {G}$, where $i_{1} < i_{2} < \ldots < i_{p}$. Let 
next$(q_i) = a_{i_1}$, next$(a_{i_p}) = q_j$ and next$(a_{i_j}) = a_{i_{j+1}}$,
for $1 \le j < p$.
Then to prove that
$X$ is an independent set in $G$, it is sufficient to show that for
each $u \in X$, $u$ is 
not anticlockwise adjacent to next$(u)$. This in fact is an easy consequence 
of Corollary~\ref {cor:fr2}. 
\end {note}

~~~~

\begin{enumerate}
\item{t is even.}


\begin{enumerate}[(i)]
\item $f(q_i)=i$ for $1 \leq i \leq r$
\item $f(a_i)= (i-1) \bmod 2x + 1$ for $1 \leq i \leq k-x+1$
\item $f(a_{k-x+i+1})=r+i$ for $1 \leq i \leq x-1$
\end{enumerate}
In order to show that $f$ is a valid coloring, we will show for each $h$,
$1\leq h \leq r+x-1$, the subset of vertices that are given color $h$  induces
an independent set. The number of 
vertices that get color $h$ vary with the different ranges of $h$ as seen below.
\begin{enumerate}
\item For the range $1 \leq h \leq x+b+1$, the arcs that get color $h$ are 
$\{ q_h, a_h, a_{h+2x},\- \ldots, a_{k-b-2x+h} \}$. Note that $((k-x+1)-1)
(\bmod 2x)+1=(k+x)(\bmod 2x)+1=x+b+1$ and hence $a_{k-x+1}$ gets color $x+b+1$. 
If $a_h$ is clockwise adjacent to $q_h$ then by Corollary~\ref{cor:fr2}.1, 
at least $r$ arcs are clockwise adjacent to $q_h$ which contradicts Lemma~\ref{lem:lem24}. Hence, $a_h$ is not clockwise adjacent to $q_h$. 
Also $a_{j+2x}$ is not clockwise adjacent to $a_j$ for $1 \leq j \leq k-3x+1$ 
as $a_{j+2x-1}$ is not clockwise adjacent to $a_j$ by Lemma~\ref{lem:lem38}. 
Arc $a_{k-b-2x+h}$ is not anticlockwise adjacent to 
$q_h$ because otherwise we have by Corollary~\ref{cor:fr2}.4, $k+1-(k-b-2x+h)+
h-1= 2x+b >r$ (since by Lemma~\ref{lem:lem39}, $2x+b>r$) arcs anticlockwise adjacent 
to $q_h$ which contradicts Lemma~\ref{lem:lem24}.
 Hence, by Note \ref {note:indset} vertices  that are given color $h$
form  an independent set.

\item For $x+b+2 \leq h \leq 2x$, the arcs that 
get color $h$ are $\{ q_h, a_h, a_{h+2x},\- \ldots, a_{k-b-4x+h} \}$. The arcs
in this set are also not adjacent to each other for the same reasons as 
discussed above.

\item For $2x < h \leq r$, only arc $q_h$ gets color $h$ and hence it is
an independent set.

\item Only one arc gets color $h$ for $h$ in the range $r < h \leq r+x-1$ 
and hence each is an independent set.
\end{enumerate}
Thus, for the case when $t$ is even,
we have demonstrated a valid vertex coloring using $r+x-1$ colors.\\


\item{t is odd.}\\
We have two sub-cases here based on the values of $x$ and $b$. We know from
Lemma~\ref{lem:lem2x} that $2x\leq r$. Also we have
$b\leq x-1$. Therefore, we have $x+b \leq r-1$ with the equality holding
when $2x=r$ and $b=x-1$.

\textbf{Sub-case 1: $x+b < r-1$} 

\begin{enumerate}[(i)]
\item $f(q_i)=i$ for $1 \leq i \leq r$
\item $f(a_i)= (i-1) (\bmod 2x) +1$ for $1 \leq i \leq k-3x-b$
\item $f(a_i)=(i-(k-2x-b))(\bmod r)+1$ for $k-2x-b \leq i \leq k-x+1$ 
\item $f(a_{k-3x-b+j})=f(a_{k-x+j+1})=r+j$ for $1 \leq j \leq x-1$
\end{enumerate}
For $1 \le h \le r+x-1$, let $X_h$ denote the set of vertices 
of $\mathcal {G}$ that are given the color $h$.  We will show 
that $X_h$ is an independent set.
As before we consider the different ranges of $h$ and study $X_h$. 
The reader  may find it  useful to note now itself that 
no two  arcs that belong  to  the range $k-2x-b \le i \le k-x+1$ 
(i.e. the range defined in  (iii) ) are given the same color,
since there are  only  at most   $(k-x+1)+1-(k-2x-b)=x+b+2 \leq r$(last inequality true as $x+b<r-1$) arcs that belong to this range. 

\begin{enumerate}
\item  For $1 \leq h \leq 2x$:  Recall that $k \ge 4x -1$. If
 $k = 4x-1 = 3x + (x-1)$ and $b=x-1$, then the range defined by
 $1 \le i \le k - 3x-b$ is empty. Thus in this case, 
 $X_h = \{q_h, a_{k-2x-b+h-1}\}$. Otherwise we consider
 two cases: If $k-2x-b+h-1 > k-x+1$ then $X_h = \{q_h,a_h, a_{h+2x},
\ldots,a_{k-5x-b+h} \}$ else
 $X_h = \{ q_h, a_h, a_{h+2x}, \ldots, a_{k-5x-b+h},a_{k-2x-b+h-1}\}$.


Now we verify that $X_h$ is an independent set in $\mathcal {G}$: 
$a_h$ is not 
clockwise adjacent to $q_h$ and $a_{j+2x}$ is not clockwise adjacent to $a_j$ 
for $1 \leq j \leq 3x$ for reasons discussed in 1(a) of $t$ is even case above.
$a_{k-5x-b+h}$ is not anticlockwise adjacent to $a_{k-2x-b+h-1}$ since the 
number of arcs anticlockwise adjacent to $a_{k-2x-b+h-1}$ would then be
$(k-2x-b+h-1)-(k-5x-b+h)=3x-1 \geq 2x+b > r$ which contradicts 
Lemma~\ref{lem:lem24}. Arc $a_{k-2x-b+h-1}$ is not anticlockwise adjacent to 
$q_h$ as the number of arcs that are anticlockwise adjacent to $q_h$ by 
Corollary~\ref{cor:fr2}.4 would then be $k+1-(k-2x-b+h-1)+(h-1) =2x+b+1 >r$ 
(since by Lemma~\ref{lem:lem39} $2x+b >r$) which contradicts 
Lemma~\ref{lem:lem24}.  Arc $q_h$ is also not anticlockwise adjacent to 
$a_{k-2x-b+h-1}$ because then by Corollary~\ref{cor:fr2}.1 number of arcs that
will be clockwise adjacent to $q_h$ would be $(k-2x-b+h-1)+(r-h)=k+(r-2x-b)-1
\geq k-(b+1) \geq k-x \geq r+x-1 \geq r$(using Lemma~\ref{lem:lem2x} and
Lemma~\ref{lem:lemr2x}) which contradicts Lemma~\ref{lem:lem24}.
Now by Note \ref {note:indset} we can easily deduce that the above set is
an independent set. 

\item For $2x < h \leq r$, arcs that get color $h$ are $q_h$ and 
$a_{k-2x-b+h-1}$ (if $k-2x-b+h-1 \leq k-x+1$). $a_{k-2x-b+h-1}$ is neither
anticlockwise adjacent nor clockwise adjacent to $q_h$ as discussed above.
Therefore, arcs that get color $h$ form an independent set.

\item Arcs that get color $r+j$ for $1\leq j\leq x-1$ are $\{a_{k-3x-b+j},
a_{k-x+j+1}\}$. Arc $a_{k-3x-b+j}$ is not anticlockwise adjacent to 
$a_{k-x+j+1}$ as then the number of arcs that are anticlockwise adjacent to 
$a_{k-x+j+1}$ by Corollary~\ref{cor:fr2}.3 would be $(k-x+j+1)-(k-3x-b+j)=
2x+b+1 >r$ which is a contradiction
by Lemma~\ref{lem:lem24}. It is easy to see that $a_{k-3x-b+j}$ is also not 
clockwise adjacent to $a_{k-x+j+1}$.
\end{enumerate}
We have seen that the arcs that get color $h$ for $1\leq h\leq r+x-1$ form
an independent set and hence the coloring is a valid coloring.

\textbf{Sub-case 2: $x+b=r-1$}

As mentioned before, if  $x+b=r-1$ we should have $2x=r$ and $b=x-1$. 
We first make the following claim:\\
\noindent \textbf{Claim}
If $x+b = r-1$, then for any $j$ in the range  $1\leq j \leq x$, 
 the arc $a_{k-2x-b+j}$  is not anticlockwise
adjacent to $a_{k-x+j}$.\\
\textbf {proof of claim:}  Suppose there exist a $j$, where
$1 \le j \le x$ such that $a_{k-2x-b+j}$ is anticlockwise adjacent to
$a_{k-x+j}$. Then, we  demonstrate a good path set, contradicting 
Lemma \ref {lem:lem37}. For $1 \le j \le x$, the $x$ sequences  defined by  
$P_j = (a_j, a_{x+j}, \ldots, a_{k-2x-b+j})$ indeed form vertex
disjoint paths (by Lemma \ref {lem:lem35}). 
We would show how to extend each path $P_j$ ($1 \le j \le x$) to 
a path $P'_{j}$ such that $P'_{j}$ ends at the vertex
$a_{k-x+j}$. If we assume that the  vertex $a_{k-2x - b +j}$
is anticlockwise adjacent to $a_{k-x+j}$ then the path $P_j$ can be
readily extended to a path $P'_{j}$, by adjoining the vertex $a_{k-x+j}$
at the end of $P_j$, i.e. just after $a_{k-2x-b+j}$. Now we extend the
remaining $x-1$ paths $P_i$  ($1\le i \le x$ and $i \ne j$) to get
$P'_{i}$ as follows:
 \begin {eqnarray}
 \mbox {~for $i < j$~~~~~~}  P'_{i} &=& (P_i,a_{k-x-b+i},a_{k-x+i}) \nonumber \\
 \mbox {~for $i > j$~~~~~~}  P'_{i} &=& (P_i,a_{k-x-b+i-1},a_{k-x+i}) \nonumber
 \end {eqnarray}
By applying Lemma \ref {lem:lem35}, it can be easily verified that 
$P'_{1},\cdots,P'_{x}$ form paths. Moreover they are vertex disjoint
paths: Notice that the last but one vertex $a_l$ in each path $P'_{i}$
($1 \le i \le x$, $i \ne j$)
belong to the range $k-x-b + 1 \le l \le k-x$. Since $b = x -1$, there
 are sufficient number of ``intermediate'' vertices to connect the last
vertex of $P_i$ (i.e. $a_{k-2x-b+i}$)   to $a_{k-x+i}$.
Thus we have demonstrated a good path set.
The claim follows.


Now, we demonstrate a vertex coloring of $\mathcal {G}$ using 
$r+x-1$ colors: 

\begin{enumerate}[(i)]
\item $f(q_i)=i$ for $1 \leq i \leq r$
\item $f(a_i)=(i-1) (\bmod 2x) +1$ for $1 \leq i \leq k-x-b$
\item $f(a_{k-x+j})=x+j$ for $1\leq j \leq x$
\item $f(a_{k-x-b+j})=r+j$ for $1 \leq j \leq b$
\end{enumerate}
As before we examine the set $X_h$, the set of arcs that get color $h$. For 
$1\leq h \leq x$, $ X_h = \{q_h, a_h, a_{h+2x},\ldots, a_{k-3x-b+h}\}$. This is 
an independent set for reasons discussed in 1(a) of $t$ is even case above.  
For $x < h \leq 2x$, $X_h = \{q_h, a_h, a_{h+2x}, \ldots, a_{k-3x-b+h}, 
a_{k-2x+h}\}$. Letting $j=h-x$, we have
$a_{k-3x-b+h} = a_{k-2x-b+j}$ which is not anticlockwise adjacent to 
$a_{k-2x+h}=a_{k-x+j}$ by the claim we proved above. 
 Also, $a_{k-2x+h}$ is not anticlockwise adjacent to $q_h$ because
otherwise the number of arcs that are anticlockwise adjacent to $q_h$ would be
$k+1-(k-2x+h)+h-1 = 2x = r$ by Corollary~\ref{cor:fr2}.4 which contradicts 
Lemma~\ref{lem:lem24}. Thus $X_h$ is an independent set for $x < h \leq 2x$.
Finally,  it is easy to see that $X_h$ is a singleton set
for  $r=2x<h\leq r+x-1$ and hence  forms an independent set.
Therefore the above coloring is a valid coloring that  
uses only $r+x-1$ colors.
\end{enumerate}

Therefore, we can see that irrespective of whether $t$ is even or odd, we can
show that we either have a good path set(and thus an $r+x$ clique minor) or we 
can color using $r+x-1$ colors which is a contradiction. Hence Hadwiger's 
conjecture is true for proper circular arc graphs.\qed
\end{pf}


\end{document}